\documentclass[a4paper,10pt]{amsart}

\usepackage{amssymb}
\usepackage{amsmath}
\DeclareMathOperator{\im}{Im}  %%Imagen n�cleo y Ext

\newcommand{\fin}{\hspace{\stretch{1}}$\square$}

\newenvironment{demostracion}{\noindent \textbf{Proof.}}

\newtheorem{definicion}{Definition}[section]

\newtheorem{teorema}[definicion]{Theorem}

\newtheorem{corolario}[definicion]{Corollary}
\newtheorem{remark}[definicion]{Remark}

\title{A note on the construction of finitely injective modules}

\author{Pedro A. Guil Asensio}
\address{Department of Mathematics. University of Murcia. 30100 Espinardo, Murcia, Spain.}
\email{paguil@um.es}
\thanks{First author has been partially supported by the DGI and by the Fundaci\'on S\'eneca. Part of the sources of both institutions come from the FEDER funds of the European Union.}

\author{Manuel C. Izurdiaga}
\address{Department of Algebra and Analysis, University of Almeria, E-04071, Almeria, Spain}
\email{mizurdia@ual.es}

\author{Blas Torrecillas}
\address{Department of Algebra and Analysis, University of Almeria, E-04071, Almeria, Spain}
\email{btorreci@ual.es}

\begin{document}

\begin{abstract}
  We develop a technique to construct finitely injective modules which
  are non trivial, in the sense that they are not direct sums of
  injective modules. As a consequence, we prove that a ring $R$ is
  left noetherian if and only if each finitely injective left
  $R$-module is trivial, thus answering an open question posed by
  Salce.
\end{abstract}

\maketitle

\section{Introduction}
\label{sec:introduction}

It has been recently shown by Salce in \cite{s} that if $R$ is a non
noetherian Matlis valuation domain, then there exists a nontrivial
finitely injective $R$-module. Where an $R$-module $M$ is called
finitely injective if each finite subset of $M$ is contained in an
injective submodule (which is necessarily a direct summand). And a
finitely injective module is said to be trivial if it is a direct sum
of injective modules. This result was inspired by an older
characterization of noetherian rings obtained by Ramamurthi and
Rangaswamy in \cite{rr}. Namely, they proved that a ring $R$ is left
noetherian if and only if every finitely injective left $R$-module is
injective. At this point, the natural question of whether general left
noetherian rings can be characterised in terms of the existence of
nontrivial finitely injective modules naturaly arises (see
\cite[Question 1]{s}).

The goal of this note is to give a positive answer to the above
question, as well as provide a simpler and more natural tool to
construct nontrivial finite injective modules over non noetherian
rings. Using this new construction, we prove that an injective left
$R$-module $M$ is $\Sigma$-injective if and only if any finitely
injective submodule of the injective envelope of $M^{(\aleph_0)}$ is
trivial. In particular, we deduce that a ring is left noetherian if
and only if any finitely injective left module is a direct sum of
injective modules.
 
Salce's construction of non-trivial finitely injective modules is
based on a classical construction by Hill for Abelian groups (see
\cite{h}). This construction was later generalized by Griffith
\cite{g} and Huisgen-Zimmermann \cite{z} for modules in order to
characterize left perfect rings in terms of the existence of
nontrivial flat and strict Mittag-Leffler (equivalently, locally
projective) modules. We would like to stress that Hill's construction
needs to be applied to countably generated modules (usually,
projective modules). Indeed, Salce's example of a non-trivial finitely
injective module is obtained by applying this Hill's construction to a
countably generated injective module which is not
$\Sigma$-injective. The existence of these countably generated modules
is guaranteed over non noetherian Matlis valuation domains, but not
over arbitrary rings.  Therefore, we need to develop in this note a
new (and much simpler) explicit construction of non trivial finitely
injective modules that can be applied to any injective module which is
not $\Sigma$-injective.

\section{Finitely injective modules over non-noetherian rings}
\label{sec:finit-inject-modul}

Along this note, $R$ will denote a ring with unit and module will mean
a left $R$-module unless otherwise is stated. Morphisms will operate
on the right and the composition of $f:A \rightarrow B$ and $g:B
\rightarrow C$ will be denoted by $fg$. Given a left $R$-module $M$,
we will denote by $E(M)$ its injective hull.

Recall that an injective module $M$ is said to be $\Sigma$-injective
if any direct sum of copies of $M$ is injective. Our main theorem
characterizes $\Sigma$-injective modules in terms of finitely
injective submodules of the injective hull of a countably direct sum
of copies of the module.

% \begin{lema} \label{l:directsummand} If $N$ is a s-pure submodule of
%   a finitely injective module then it is finitely injective.
% \end{lema}

\begin{teorema}\label{t:maintheorem}
  Let $M$ be an injective module. The following assertions are
  equivalent:
  \begin{enumerate}
  \item $M$ is $\Sigma$-injective.

  \item Every finitely injective submodule of
    $E\left(M^{(\aleph_0)}\right)$ is injective.

  \item Every finitely injective submodule of
    $E\left(M^{(\aleph_0)}\right)$ is a direct sum of injective
    modules.
  \end{enumerate}
\end{teorema}

\begin{demostracion}
  $1)\Rightarrow 2)$. Let us first note that
  $E\left(M^{(\aleph_0)}\right) = M^{(\aleph_0)}$ is
  $\Sigma$-injective too and, consequently, it is a direct sum of
  indecomposable modules (see \cite{c}). Now, using \cite[Theorem
  2.22]{mm} and \cite[Lemma 2.16]{mm} we get that the union of any
  chain of direct summands of $M^{(\aleph_0)}$ is a direct summand.

  Let $N$ be a finitely injective submodule of $M^{(\aleph_0)}$ and
  $\{x_\alpha\,|\,\alpha < \kappa\}$ be any generating set of $N$
  (where $\kappa$ is an ordinal). We claim that $N$ is the union of a
  chain $\{N_\alpha\,|\,\alpha < \kappa\}$ of direct summands of
  $M^{(\aleph_0)}$. By the previous discussion, this implies that $N$
  is a direct summand of $M^{(\aleph_0)}$ and, in particular,
  injective. We are going to construct this chain of submodules of $N$
  recursively on $\alpha$ with the property that $x_\alpha \in
  N_\alpha$ for each $\alpha < \kappa$.
  
  For $\alpha=0$, choose an injective submodule $N_0$ of $N$ that
  contains $x_0$. We know that this $N_0$ does exist since we are
  assuming that $N$ is finitely injective. As $N_0$ is injective, it
  is a direct summand of $M^{(\aleph_0)}$ and therefore, of $N$.
  
  Let now $\alpha<\kappa$ be any ordinal and assume that we have
  constructed our chain $\{N_\gamma\,|\,\gamma < \alpha\}$ of
  submodules of $N$. Then note that $\bigcup_{\gamma <
    \alpha}N_\gamma$ is a direct summand of $M^{(\aleph_0)}$ (and
  therefore, of $N$), since it is the union of a chain of direct
  summands of $M$ and $M$ is $\Sigma$-injective. Let us write $N =
  \left(\bigcup_{\gamma < \alpha}F_\gamma\right) \oplus E$ for some $E
  \leq F$, and $x_\alpha = n + e$ for some $n \in \bigcup_{\gamma <
    \alpha}N_\gamma$ and $e \in E$. Since $E$ is finitely injective,
  there exists an injective submodule $L$ of $E$ containing $e$. Let
  $E(Re)$ be an injective hull of $Re$ inside $L$. Then $E(Re)\cap
  (\bigcup_{\gamma < \alpha}N_\gamma)=0$ since $Re$ is essential in
  $E(Re)$ and $Re\cap (\bigcup_{\gamma < \alpha}N_\gamma)=0$. Set then
  $N_\alpha = \left(\cup_{\gamma < \alpha}N_\gamma\right) \oplus
  E(Re)$ which contains $x_\alpha$ and, by \cite[Proposition 2.2]{mm},
  it is a direct summand of $M^{(\aleph_0)}$. This concludes the
  construction.

  $2) \Rightarrow 3)$. Trivial.

  $3) \Rightarrow 1)$. Suppose that $(1)$ is false. We are going to
  construct a finitely injective submodule of $E(M^{(\aleph_0)})$
  which is not a direct sum of injective modules. As we are assuming
  that $M$ is not $\Sigma$-injective, we know that $M^{(\aleph_0)}$ is
  not injective by \cite[Proposition 3]{f}. Denote by
  $N=M^{(\aleph_0)}$ and by $E=E\left(M^{(\aleph_0)}\right)$. By
  Baer's criterium, there exists a left ideal $I$ of $R$ and a
  morphism $f:I \rightarrow N$ that cannot be extended $R$. Let
  $\Omega$ be the set of all submodules $L$ of $E$ such that:
  \begin{enumerate}
  \item $L$ is finitely injective;

  \item $N \leq L$, and

  \item The morphism $f:I \rightarrow L$ cannot be extended to a
    morphism $R\rightarrow L$.
  \end{enumerate}

  Clearly, $\Omega$ is a non-empty partially ordered set. Let us show
  that it is inductive. Let $\{L_k:k \in K\}$ be a chain in $\Omega$
  and let us prove that $L=\bigcup_{k \in K}L_k \in
  \Omega$. Trivially, $L$ satisfies conditions $(1)$ and $(2)$. Suppose
  that there exists an extension $g:R \rightarrow L$ of $f$. Then, as
  $R$ is finitely generated, there exists a $k \in K$ such that $\im f
  \leq L_k$. But this means that $g$ is an extension of $f$ over
  $L_k$, which contradicts that $L_k \in \Omega$. Thus, $L$ also
  satisfies $(3)$ and it is an element of $\Omega$.

  Let now $L$ be a maximal element of $\Omega$. We are going to show
  that $L$ cannot be a direct sum of injective modules. Suppose on the
  contrary that $L$ would be a direct sum of injective modules, say $L
  = \oplus_{t \in T}E_t$. Denote by $q_t:L \rightarrow E_t$ the
  canonical projection, for each $t \in T$. Let
  \begin{displaymath}
    T' = \{t \in T: fq_t \neq 0\}
  \end{displaymath}
  and note that $T'$ is infinite because otherwise, $\im f$ would be
  contained in a finite direct subsum of $\oplus_{t \in T}E_t$ and $f$
  would have an extension to $R$. Call $q_{T'}:L \rightarrow \oplus_{t
    \in T'}E_{t}$ the projection. And write $T' = T_1 \cup T_2$ as a
  disjoint union of two infinite subsets. Denote by
  $q_{T_i}:\oplus_{t \in T'}E_t \rightarrow \oplus_{t \in T_i}E_t$ the
  projections, for $i = 1,2$.

  We claim that neither $f q_{T'}q_{T_1}$ nor $f q_{T'}q_{T_2}$ can be
  extended to $R$. Assume on the contrary that, for instance, $f
  q_{T'} q_{T_1}$ could be extended to a morphism $h:R \rightarrow
  \oplus_{t \in T_1}E_t$. Then there would exist a finite subset $F
  \subseteq T_1$ such that $\im h \subseteq \oplus_{t \in F}E_t$,
  as $\im h$ is finitely generated. But then $0 = f q_{T'} q_{T_1}
  q_t = f q_t$ for every $t \in T_1 \setminus F$. A contradiction that
  proves our claim.

  Let finally $L' = E(\oplus_{t \in T_1}E_t)\bigoplus\left(\oplus_{t
      \in T \setminus T_1}E_t\right)$. Then $L'$ belongs to $\Omega$,
  since $f$ does not have an extension to $R$ because
  $T_2 \subseteq T \setminus T_1$ and $f q_{T'} q_{T_2}$ cannot be
  extended to $R$. Moreover, $L$ is strictly contained in $L'$, as
  $\oplus_{t \in T_1}E_t$ is not injective. But this contradicts the
  maximality of $L$.\fin
\end{demostracion}

\begin{remark}
  Let us note that our proof of $(3) \Rightarrow (1)$ in the above
  theorem gives a general tool to construct finitely injective modules
  which cannot be injective. This construction is simpler than the one
  obtained \cite{s} and seems more natural on this injectivity
  context.
\end{remark}

As a consequence, we get a possitive answer to \cite[Question
1]{s}.

\begin{corolario}
  Let $R$ be a ring. The following assertions are equivalent:
  \begin{enumerate}
  \item $R$ is left noetherian.

  \item Each finitely injective left $R$-module is injective.

  \item Each finitely injective left $R$-module is a direct sum of
    injective modules.
  \end{enumerate}
\end{corolario}


\begin{thebibliography}{99}
\bibitem{c} A. Cailleau: Une caractérisation des modules
  $\Sigma$-injectifs, C. R. Acad. Sci. Paris Ser A 269 (1969),
  997-999.

\bibitem{f} C. Faith. Rings with ascending chain conditions on
  annihilators. Nagoya Math. J. 27 (1966), 179-191.

\bibitem{g} P. Griffith. A note on a theorem of Hill, Pacific
  J. Math. 29 (1969), 279-284.

\bibitem{h} P. Hill. On the decomposition of groups,
  Canad. J. Math. 21 (1969), 762-768.

\bibitem{mm} S. H. Mohamed, Bruno J. M\"uller, Continuous and discrete
  modules.  London Mathematical Society Lecture Note Series,
  147. Cambridge University Press, Cambridge, 1990. viii+126 pp. ISBN:
  0-521-39975-0.

\bibitem{rr} V. S. Ramamurthi, K. M. Rangaswamy. On finitely injective
  modules, J. Austral. Math. Soc. 16 (1973), 239-248.

\bibitem{s} Luigi Salce. On Finitely Injective Modules and Locally
  Pure-Injective Modules over Prüfer
  Domanis. Proc. Amer. Math. Soc. Vol. 135, n0 11, November 2007,
  Pages 3485-3493.

\bibitem{z} B. Zimmermann-Huisgen. On the abundance of
  $\aleph_1$-separable modules, Abelian groups and noncommutative
  rings, Contemp. Math., vol. 130, Amer. Math. Soc., Provicence, RI,
  1992, pp. 167-180.
\end{thebibliography}
\end{document}